\documentclass[psamsfonts,reqno]{amsart}
\usepackage{amssymb,eucal,latexsym,xy,graphics}
\xyoption{all}
\xyoption{ps}

\theoremstyle{plain}

\newtheorem{lemma}{Lemma}[section]
\newtheorem{theorem}[lemma]{Theorem}
\newtheorem{proposition}[lemma]{Proposition}

\newtheorem*{sclaim}{Claim}

\newtheorem*{stat}{\name}
\newcommand{\name}{testing}

\theoremstyle{definition}
\newtheorem{definition}[lemma]{Definition}
\newtheorem{problem}{Problem}
\newtheorem{example}[lemma]{Example}

\theoremstyle{remark}

\newcommand{\qedc}{{\qed}~{\rm Claim~{\theclaim}.}}
\newcommand{\qedsc}{{\qed}~{\rm Claim.}}

\newenvironment{scproof}
{\begin{proof}[Proof of Claim.]}
{\qedsc\renewcommand{\qed}{}\end{proof}}

\newcommand{\case}[1]%
{\smallskip\noindent\textbf{\textit{Case}\ {#1}.}}

\numberwithin{equation}{section}
\numberwithin{figure}{section}

\newcommand{\pup}[1]{\textup{(}{#1}\textup{)}}

\DeclareMathOperator{\Retr}{Retr}
\DeclareMathOperator{\Ob}{Ob}

\DeclareMathOperator{\Idc}{Id_c}

\DeclareMathOperator{\Conc}{Con_c}

\DeclareMathOperator{\FP}{FP}
\DeclareMathOperator{\Dim}{Dim}

\newcommand{\la}{\boldsymbol{a}}

\newcommand{\lx}{\boldsymbol{x}}
\newcommand{\ly}{\boldsymbol{y}}

\newcommand{\lu}{\boldsymbol{u}}

\newcommand{\bL}{\boldsymbol{L}}

\newcommand{\xD}{\mathbf{D}}
\newcommand{\xE}{\mathbf{E}}
\newcommand{\xF}{\mathbf{F}}
\newcommand{\xQ}{\mathbf{Q}}
\newcommand{\xR}{\mathbf{R}}

\newcommand{\cA}{\mathcal{A}}
\newcommand{\cB}{\mathcal{B}}
\newcommand{\cC}{\mathcal{C}}
\newcommand{\cD}{\mathcal{D}}

\newcommand{\cG}{\mathcal{G}}
\newcommand{\cM}{\mathcal{M}}
\newcommand{\cR}{\mathcal{R}}
\newcommand{\cS}{\mathcal{S}}
\newcommand{\cI}{\mathcal{I}}

\newcommand{\cV}{\mathcal{V}}

\newcommand{\QQ}{\mathbb{Q}}

\newcommand{\bbD}{\mathbb{D}}
\newcommand{\tbbD}{\widetilde{\mathbb{D}}}
\newcommand{\hxD}{\widehat{\xD}}
\newcommand{\bbE}{\mathbb{E}}

\newcommand{\cm}{commutative monoid}

\newcommand{\eps}{\varepsilon}

\newcommand{\into}{\hookrightarrow}
\newcommand{\onto}{\twoheadrightarrow}

\newcommand{\seq}[1]{\left\langle{#1}\right\rangle}
\newcommand{\seqm}[2]{\seq{#1\mid#2}}

\newcommand{\set}[1]{\{#1\}}
\newcommand{\setm}[2]{\set{#1\mid#2}}

\newcommand{\go}{\omega}
\newcommand{\gos}{\omega\setminus\set{0}}

\newcommand{\ol}[1]{\overline{#1}}
\newcommand{\oll}[1]{\,\overline{\!#1}}

\newcommand{\id}{\mathrm{id}}
\newcommand{\jz}{$\langle\vee,0\rangle$}
\newcommand{\jzu}{$\langle\vee,0,1\rangle$}
\newcommand{\bl}{$\langle\vee,\wedge,0,1\rangle$}
\newcommand{\jzs}{\jz-semi\-lat\-tice}
\newcommand{\jzus}{\jzu-semi\-lat\-tice}
\newcommand{\jzh}{\jz-ho\-mo\-mor\-phism}

\newcommand{\jze}{\jz-em\-bed\-ding}
\newcommand{\jzue}{\jzu-em\-bed\-ding}
\newcommand{\ble}{\bl-em\-bed\-ding}

\newcommand{\fin}{\mathrm{fin}}
\newcommand{\vpd}{\mathbin{\Pi}}

\begin{document}

\title[Lifting retractable diagrams]%
{Lifting retracted diagrams with respect to projectable
functors}

\author[F.~Wehrung]{Friedrich Wehrung}
\address{LMNO, CNRS UMR 6139\\ D\'epartement de Math\'ematiques\\
Universit\'e de Caen\\ 14032 Caen Cedex\\ France}
\email{wehrung@math.unicaen.fr}

\urladdr{http://www.math.unicaen.fr/\~{}wehrung}

\date{\today}



\thanks{The author was partially supported by the institutional
grant CEZ:J13/98:1132000007a, by FRVS 2731/2003, by the Fund
of Mobility of the Charles University (Prague), and by INTAS
project 03-51-4110}

\begin{abstract}
We prove a general categorical theorem that enables us to state
that under certain conditions, the range of a functor is large. As
an application, we prove various results of
which the following is a prototype: \emph{If
every diagram, indexed by a lattice, of
finite \emph{Boolean} \jzs s with \jze s, can be lifted with
respect to the $\Conc$ functor on lattices, then so can every
diagram, indexed by a lattice, of finite
\emph{distributive} \jzs s with \jze s}. If the premise of this
statement held, this would solve in turn the (still open) problem
whether every distributive algebraic lattice is isomorphic to the
congruence lattice of a lattice. We also outline potential
applications of the method to other functors, such as the
$R\mapsto V(R)$ functor on von~Neumann regular rings.
\end{abstract}

\maketitle

\section{Introduction}\label{S:Intro}

We are dealing with the following general kind of problem.
We are given a functor $\mathbf{F}$ from a category $\cA$ to a
category $\cB$, we wish to investigate the range of $\xF$, that
is, the class of all objects of $\cB$ isomorphic to some $\xF(A)$,
where $A$ is an object of~$\cA$. The goal of the present paper is
to provide a general categorical framework for proving results
saying that if certain diagrams of $\cB$ can be lifted with
respect to~$\xF$, then the range of $\xF$ is large.

Let us be more precise. We assume that $\cB$ has a subcategory
$\cB'$, in which many diagrams can be lifted with respect to
$\xF$. Furthermore, we assume that every member of $\cB$ is a
\emph{retract} of some member of $\cB'$, and this in a functorial
way; we say that $\cB$ is a \emph{functorial retract} of $\cB'$
(see Definition~\ref{D:FunctRetr}). Then, under certain
conditions, many more diagrams in $\cB$ (not only in $\cB'$)
can be lifted with respect to~$\xF$. This is the main result of
the paper, Theorem~\ref{T:Main}. It is stated in the language of
category theory.

We do not challenge the view that Theorem~\ref{T:Main} is a
technical, at first sight unappealing result,
whose formulation lies apparently far away from most
readers' usual mathematical practice. However, the main motivation
for establishing it lies in an easily formulated problem,
namely R.\,P. Dilworth's Congruence Lattice Problem CLP, open
since 1945, that asks whether any distributive algebraic lattice
is isomorphic to the congruence lattice of a lattice, see
G. Gr{\"a}tzer and E.\,T. Schmidt~\cite{GrScC} or
J. T\r{u}ma and F. Wehrung~\cite{CLPSurv} for surveys.
Theorem~\ref{T:Main} makes it possible to reduce a
certain diagrammatic strengthening of CLP to diagrams of
\emph{finite Boolean} semilattices and \emph{\jze s}, see
Proposition~\ref{P:UltrabRed}. This reduction requires the result,
established by the author in
\cite{Ultra}, that states that \emph{every finite distributive
\jzs\ is a \jz-retract, in a functorial way, of some finite Boolean
\jzs}. Moreover, the categorical formulation of
Theorem~\ref{T:Main} makes it possible to attack other
representation problems from a new perspective: for
example, functorial lifting, with respect to the compact
congruence semilattice functor $\Conc$, of finite \jzs s by
algebras of a given similarity type (see Theorem~\ref{T:LiftConVar}
and Proposition~\ref{P:NonLiftSig}), or representation problem of
countable refinement monoids by the nonstable K-theory of
von~Neumann regular rings (see Proposition~\ref{P:V(R)attempt}).
Although no definite answer of the corresponding long-standing
open problems is reached here, we hope that the ideas of the
present paper will provide a new start towards more final
conclusions.

We conclude the paper by mentioning further functors for which our
results might be of some potential use, together with a few open
problems.

While the present paper deals with how to \emph{use} functorial
retractions in order to prove that certain functors have large
range, the paper~\cite{Ultra} deals with the \emph{existence}
of functorial retractions.

As for \cite{Ultra}, the techniques introduced in
the present paper were developed to solve problems
closer to universal algebra than to category theory. Nevertheless,
the paper is, up to Section~\ref{S:LiftxD}, purely
category-theoretical. For these reasons, the author chose to write
it in probably more detail than a category theorist would wish,
with the hope to make it reasonably intelligible to members of both
communities.

\section{Basic concepts}\label{S:Basic}

We shall identify every natural number $n$ with the set
$\set{0,1,\dots,n-1}$. Similarly,
$\go=\set{0,1,2,\dots}$, the set of all natural numbers.

We shall mostly use the notation and terminology from S. Mac
Lane~\cite{McLa} and from the author's previous paper
\cite{Ultra}. A morphism in a given category is an \emph{epic}, if
it is right cancellable; we shall often denote epics by double
arrows, as in $f\colon A\onto B$.

The following definitions from \cite{Ultra} will be crucial.

\begin{definition}\label{D:RetrAB}
Let $\cA$ and $\cB$ be
subcategories of a category $\cC$. We denote by $\Retr(\cA,\cB)$
the category whose objects and morphisms are the following:
\begin{itemize}
\item[---] \emph{Objects}: all quadruples $\seq{A,B,\eps,\mu}$,
where
$A\in\Ob\cA$, $B\in\Ob\cB$, $\eps\colon A\to\nobreak B$,
$\mu\colon B\to A$, and $\mu\circ\eps=\id_A$.

\item[---] \emph{Morphisms}: a morphism from $\seq{A,B,\eps,\mu}$
to $\seq{A',B',\eps',\mu'}$ is a pair $\seq{f,g}$, where
$f\colon A\to A'$ in $\cA$, $g\colon B\to B'$ in $\cB$,
$g\circ\eps=\eps'\circ f$, and $\mu'\circ g=f\circ\mu$.
Composition of morphisms is defined by the rule
$\seq{f',g'}\circ\seq{f,g}=\seq{f'\circ f,g'\circ g}$.
\end{itemize}
The \emph{projection functor} from $\Retr(\cA,\cB)$ to $\cA$ is
the functor from $\Retr(\cA,\cB)$ to~$\cA$ that sends any object
$\seq{A,B,\eps,\mu}$ to $A$ and any morphism
$\seq{f,g}$ to $f$.
\end{definition}

\begin{definition}\label{D:FunctRetr}
We say that $\cA$ is a
\emph{functorial retract} of $\cB$, if the projection functor from
$\Retr(\cA,\cB)$ to~$\cA$ has a right inverse. We shall call such
an inverse a \emph{functorial retraction} of $\cA$ to $\cB$.
\end{definition}

Hence a functorial retraction may be viewed as a triple
$\seq{\Phi,\eps,\mu}$ that satisfies the following conditions:
\begin{itemize}
\item[---] $\Phi$ is a functor from $\cA$ to $\cB$.

\item[---] For every morphism $f\colon X\to Y$ in $\cA$, we have
$\eps_X\colon X\to\Phi(X)$, $\mu_X\colon\Phi(X)\to\nobreak X$,
$\mu_X\circ\eps_X=\id_X$, $\Phi(f)\circ\eps_X=\eps_Y\circ f$, and
$\mu_Y\circ\Phi(f)=f\circ\mu_X$.
\end{itemize}

We shall often deal with colimits with respect to upward directed
posets, in short \emph{directed colimits}
(often called \emph{direct limits} in other mathematical
areas). We say that a category $\cA$ has
\emph{directed $\go$-colimits}, if every diagram of $\cA$ of the
form
 \[
 {
 \def\labelstyle{\displaystyle}
 \xymatrix{
 A_0\ar[r]^{f_0} & A_1\ar[r]^{f_1} & A_2\ar[r]^(.4){f_2} & 
 \cdots\ \cdots
 }
 }
 \]
(indexed by the poset $\go$) has a colimit in $\cA$. For a category
$\cB$, we say that a functor $\xF\colon\cA\to\cB$ \emph{preserves
directed $\go$-colimits}, if whenever $A=\varinjlim_{n<\go}A_n$,
with transition morphisms $f_n\colon A_n\to A_{n+1}$ and limiting
morphisms $g_n\colon A_n\to A$ (for $n<\go$), the relation
$\xF(A)=\varinjlim_{n<\go}\xF(A_n)$ holds, with transition
morphisms $\xF(f_n)$ and limiting morphisms $\xF(g_n)$ (for
$n<\go$).

The following easy result is in some sense the central idea that
makes it possible in this paper to express objects as
colimits of ``simpler'' objects.

\begin{lemma}\label{L:DirLimIdemp}
Let $\cB$ be a category and let $\seq{A,B,\eps,\mu}$ be an object
of $\Retr(\cB,\cB)$. Put $\rho=\eps\circ\mu$. Then $A$ is the
colimit of the sequence
 \[
 {
 \def\labelstyle{\displaystyle}
 \xymatrix{
 B\ar[r]^{\rho} & B\ar[r]^{\rho} & B\ar[r]^(.4){\rho} & 
 \cdots\ \cdots
 }
 }
 \]
with constant limiting morphism $\mu\colon B\to A$.
\end{lemma}

\begin{proof}
Clearly, $\rho\circ\rho=\rho$ and $\mu\circ\rho=\mu$. Now let $C$
be an object of $\cB$ and let $\seqm{\varphi_n}{n<\go}$ be a
sequence of morphisms, $\varphi_n\colon B\to C$, such that
$\varphi_n=\varphi_{n+1}\circ\rho$ for all $n<\go$. Since $\rho$
is idempotent, $\varphi_n=\varphi_0$ for all $n<\go$.
The morphism $\varphi_0\circ\eps$ is the unique morphism
$\psi\colon A\to C$ such that $\psi\circ\mu=\varphi_0$.
\end{proof}

We shall use the following notation for \emph{products}. For
objects $A$, $B$, and $C$ in a given category, the notation
$C=A\vpd B$ (wrt. $a$, $b$) means that $C$ is the product of
$A$ and $B$ with projections $a\colon C\to A$ and
$b\colon C\to B$. For an object $D$ and morphisms
$f\colon D\to A$ and $g\colon D\to B$, we shall denote by
$f\times g$ the unique morphism $h\colon D\to C$ such that
$f=a\circ h$ and $g=b\circ h$. In case $D=A'\vpd B'$ (wrt. $a'$,
$b'$), $f\colon A'\to A$, and $g\colon B'\to B$, we shall put
$f\vpd g=(f\circ a')\times(g\circ b')$.

If $\cA$, $\cB$, and $\cI$ are categories, a \emph{lifting} of a
functor $\xD\colon\cI\to\cB$ with respect to a functor
$\xF\colon\cA\to\cB$ is a functor $\xE\colon\cI\to\cA$ such that
the composition $\xF\xE$ is equivalent to $\xD$.

Our \cm s will be denoted additively. Let $M$ be a \cm. We
say that $M$ is a \emph{refinement monoid} (see
H. Dobbertin~\cite{Dobb83} or F. Wehrung~\cite{Wehr92a,Wehr92b}),
if $a_0+a_1=b_0+b_1$ in $M$ implies that there are elements
$c_{i,j}\in M$ (for $i$, $j<2$) such that
$a_i=c_{i,0}+c_{i,1}$ and $b_i=c_{0,i}+c_{1,i}$, for all $i<2$. We
say that $M$ is \emph{conical}, if $x+y=0$ implies that $x=y=0$,
for all $x$, $y\in M$. It is well-known that a
\jzs\ is a refinement monoid if{f} it is distributive, see
G. Gr\"atzer~\cite[Section~II.5]{GLT2}.
Every \cm\ is endowed with its \emph{algebraic quasi-ordering},
defined by $x\leq y$ $\Leftrightarrow$ $(\exists z)(x+z=y)$. A
monoid homomorphism $f\colon M\to N$ is an \emph{embedding}, if it
is one-to-one and $f(x)\leq f(y)$ if{f} $x\leq y$, for all $x$,
$y\in M$. An \emph{order-unit} of
$M$ is an element $u\in M$ such that for all $x\in M$, there
exists a positive integer $n$ such that $x\leq nu$.

A subset $I$ of a \cm\ $M$ is an \emph{o-ideal}, if $0\in I$ and
$x+y\in I$ if{f} $x\in\nobreak I$ and $y\in I$, for all $x$,
$y\in M$. For $x$, $y\in M$, let $x\equiv_Iy$ hold, if there are
$u$, $v\in I$ such that $x+u=y+v$. Then $\equiv_I$ is a monoid
congruence on $M$; we put $M/I=M/{\equiv_I}$.
A monoid homomorphism $f\colon M\to N$ between \cm s $M$ and $N$
is \emph{ideal-induced}, if $f$ is \emph{surjective} and the
o-ideal $I=f^{-1}\set{0}$ satisfies the condition that $f(x)=f(y)$
if{f} $x\equiv_Iy$, for all $x$, $y\in M$. Equivalently, the
homomorphism $f$ is isomorphic, in the category of homomorphisms
with domain~$M$, to the canonical projection $M\onto M/I$.

\section{Projecting morphisms with respect to projectable
functors}\label{S:DecDistr}

\begin{definition}\label{D:Proj}
In a category $\cC$, a morphism $f\colon A\to B$ is a
\emph{projection}, if there exists a decomposition
of the form $A=B\vpd B'$ (wrt. $f$, $f'$) in $\cC$.
\end{definition}

\begin{definition}\label{D:ProjFunct}
Let $\xF$ be a functor from a category $\cA$ to a category $\cB$,
let $A\in\Ob\cA$, $B\in\Ob\cB$, and $\varphi\colon\xF(A)\to B$.
A \emph{projectability witness} for $\seq{\varphi,A,B}$ (with
respect to~$\xF$) is a pair $\seq{a,\eps}$ satisfying the following
conditions:
\begin{enumerate}
\item $a\colon A\onto\oll{A}$ is an epic in $\cA$.

\item $\eps\colon\xF(\oll{A})\to B$ is an isomorphism in $\cB$.

\item $\varphi=\eps\circ\xF(a)$.

\item For every $f\colon A\to X$ in $\cA$ and every
$\eta\colon\xF(\oll{A})\to\xF(X)$ such that
$\xF(f)=\eta\circ\xF(a)$, there exists $g\colon\oll{A}\to X$ in
$\cA$ such that $f=g\circ a$ and $\eta=\xF(g)$.
\end{enumerate}

We say that $\xF$ is \emph{projectable}, if every triple
$\seq{\varphi,A,B}$, where $\varphi\colon\xF(A)\to B$ is a
\emph{projection}, has a projectability witness.
\end{definition}

We observe that the morphism $g$ in (iii) above is necessarily
\emph{unique} (because~$a$ is an epic). Furthermore, the
projectability witness $\seq{a,\eps}$ is unique up to isomorphism.

Definition~\ref{D:ProjFunct} is illustrated on
Figure~\ref{Fig:ProjxF}.

\begin{figure}[htb]
 \[
 {
 \def\labelstyle{\displaystyle}
 \xymatrix{
 \xF(A)\ar[r]^{\varphi}\ar[d]_{\xF(a)} & B &
 \xF(A)\ar[d]_{\xF(a)}\ar[drr]^{\xF(f)} & & &
 A\ar@{->>}[d]_a\ar[dr]^f & \\
 \xF(\oll{A})\ar[ru]|-{\eps\cong} & &
 \xF(\oll{A})\ar[rr]_{\eta=\xF(g)} & & \xF(X)
 & \oll{A}\ar[r]_g & X
 }
 }
 \]
\caption{Projectability of the functor $\xF$.}
\label{Fig:ProjxF}
\end{figure}

\section{{}From $\bbD\colon\cI\to\Retr(\cB,\cB)$ to the unfolding
$\tbbD\colon\cI\times\go\to\cB$}\label{S:D2tD}

In this section we shall fix categories $\cI$ and $\cB$, together
with a functor\linebreak $\bbD\colon\cI\to\Retr(\cB,\cB)$
(see Definition~\ref{D:RetrAB}). We also assume that $\cB$ has
nonempty finite products. (Equivalently, $B_0\vpd B_1$ exists, for
all objects $B_0$ and $B_1$ of $\cB$.)
We shall describe a construction that creates a functor
$\tbbD\colon\cI\times\go\to\cB$, that we shall call the
\emph{unfolding} of $\bbD$.

The functor~$\bbD$ is given by functors $\xD$,
$\hxD\colon\cI\to\cB$, together with correspondences
$X\mapsto\eps_X$ and $X\mapsto\mu_X$, for
$X\in\Ob\cI$, such that the following conditions hold for all
$f\colon X\to Y$ in $\cI$ (see Figure~\ref{Fig:FunctbbD}):
\begin{enumerate}
\item $\mu_X\circ\eps_X=\id_{\xD(X)}$;

\item $\hxD(f)\circ\eps_X=\eps_Y\circ\xD(f)$ and
$\mu_Y\circ\hxD(f)=\xD(f)\circ\mu_X$.
\end{enumerate}

\begin{figure}[htb]
 \[
 {
 \def\labelstyle{\displaystyle}
 \xymatrix{
 \hxD(X)\ar[r]^{\hxD(f)}\ar@{->>}[d]<.5ex>^{\mu_X} &
 \hxD(Y)\ar@{->>}[d]<.5ex>^{\mu_Y}\\
 \xD(X)\ar@<.5ex>[u]^{\eps_X}\ar[r]_{\xD(f)} &
 \xD(Y)\ar@<.5ex>[u]^{\eps_Y}
 }
 }
 \]
\caption{Description of the functor $\bbD$.}
\label{Fig:FunctbbD}
\end{figure}

In particular, $\rho_X=\eps_X\circ\mu_X$ is an idempotent
endomorphism of $\hxD(X)$. We shall need the following lemma.

\begin{lemma}\label{L:D(f)rho}
For any $f\colon X\to Y$ in $\cI$, the equality
$\hxD(f)\circ\rho_X=\rho_Y\circ\hxD(f)$ holds.
\end{lemma}

\begin{proof}
We compute $\hxD(f)\circ\eps_X\circ\mu_X=
\eps_Y\circ\xD(f)\circ\mu_X=\eps_Y\circ\mu_Y\circ\hxD(f)$.
\end{proof}

For $X\in\Ob\cI$, we define $\hxD^n(X)$ by induction on $n$, by
putting
 \begin{align}
 \hxD^1(X)&=\hxD(X)\label{Eq:hxD1X}\\
 \hxD^{n+1}(X)&=\hxD(X)\vpd\hxD^n(X)
 \quad(\text{wrt. }\alpha^X_{n+1},\,\pi^X_{n+1})\label{Eq:hxDn+1X}
 \end{align}
We also put $\alpha^X_1=\id_{\hxD(X)}$, while $\pi^X_1$ is
undefined. We also put
 \[
 \sigma^X_n=(\rho_X\circ\alpha^X_n)\times\id_{\hxD^n(X)},
 \]
with respect to the decomposition~\eqref{Eq:hxDn+1X}
(see the left half of Figure~\ref{Fig:hDn(f)sig}).

If $f\colon X\to Y$ in $\cI$,
we define $\hxD^n(f)$ by induction on $n$, by putting
$\hxD^1(f)=\hxD(f)$ and $\hxD^{n+1}(f)=\hxD(f)\vpd\hxD^n(f)$, with
respect to $\alpha^X_{n+1}$, $\pi^X_{n+1}$. Hence
$\hxD^n(f)\colon\hxD^n(X)\to\hxD^n(Y)$ in $\cB$
(see the right half of Figure~\ref{Fig:hDn(f)sig}).

\begin{figure}[htb]
 \[
 {
 \def\labelstyle{\displaystyle}
 \xymatrix{
 & \hxD^n(X)\ar[ddl]_{\rho_X\circ\alpha^X_n}\ar[d]|-{\ \sigma^X_n}
 \ar[ddr]^{\id_{\hxD^n(X)}} & & &
 \hxD^{n+1}(X)\ar[dl]_{\alpha^X_{n+1}}\ar[d]|-{\hxD^{n+1}(f)}
 \ar[dr]^{\pi^X_{n+1}} & \\
 & \hxD^{n+1}(X)
 \ar[dl]^{\alpha^X_{n+1}}\ar[dr]_{\pi^X_{n+1}} & &
 \hxD(X)\ar[d]_{\hxD(f)} &
 \hxD^{n+1}(Y)
 \ar[dl]^{\alpha^Y_{n+1}}\ar[dr]_{\pi^Y_{n+1}} &
 \hxD^n(X)\ar[d]^{\hxD^n(f)}\\
 \hxD(X) & & \hxD^n(X) & \hxD(Y) & & \hxD^n(Y)
 }
 }
 \]
\caption{The morphisms $\sigma^X_n$ and $\hxD^n(f)$.}
\label{Fig:hDn(f)sig}
\end{figure}

The following lemma is a trivial consequence of the relations
on the left half of Figure~\ref{Fig:hDn(f)sig}.

\begin{lemma}\label{L:sigoalph}
For any $X\in\Ob\cI$ and $n\in\gos$, the equality
$\alpha^X_{n+1}\circ\sigma^X_n=\rho_X\circ\alpha^X_n$ holds.
\end{lemma}

\begin{lemma}\label{L:alphD(f)}
For any $f\colon X\to Y$ in $\cI$ and $n\in\gos$, the equality
$\alpha^Y_n\circ\hxD^n(f)=\hxD(f)\circ\alpha^X_n$ holds.
\end{lemma}

\begin{proof}
For $n=1$ it is trivial. At stage $n+1$, with $n>0$, the
conclusion follows from the relations on the right half of
Figure~\ref{Fig:hDn(f)sig}.
\end{proof}

\begin{lemma}\label{L:sigD(f)}
For $f\colon X\to Y$ in $\cI$ and $n\in\gos$, the equality
$\hxD^{n+1}(f)\circ\sigma^X_n=\sigma^Y_n\circ\hxD^n(f)$ holds.
\end{lemma}

\begin{proof}
By the definition of the statement \eqref{Eq:hxDn+1X},
it suffices to prove that the compositions of the desired equality
on the left with both $\alpha^Y_{n+1}$ and $\pi^Y_{n+1}$ hold.

\noindent Composition with $\alpha^Y_{n+1}$: we calculate, using
Figure~\ref{Fig:hDn(f)sig},
 \begin{align*}
 \alpha^Y_{n+1}\circ\hxD^{n+1}(f)\circ\sigma^X_n
 &=\hxD(f)\circ\alpha^X_{n+1}\circ\sigma^X_n\\
 &=\hxD(f)\circ\rho_X\circ\alpha^X_n\\
 &=\rho_Y\circ\hxD(f)\circ\alpha^X_n
 &&(\text{by Lemma~\ref{L:D(f)rho}})\\
 &=\rho_Y\circ\alpha^Y_n\circ\hxD^n(f)
 &&(\text{by Lemma~\ref{L:alphD(f)}})\\
 &=\alpha^Y_{n+1}\circ\sigma^Y_n\circ\hxD^n(f).
 \end{align*}

\noindent Composition with $\pi^Y_{n+1}$: we calculate, using
Figure~\ref{Fig:hDn(f)sig},
 \begin{align*}
 \pi^Y_{n+1}\circ\hxD^{n+1}(f)\circ\sigma^X_n
 &=\hxD^n(f)\circ\pi^X_{n+1}\circ\sigma^X_n\\
 &=\hxD^n(f)\\
 &=\pi^Y_{n+1}\circ\sigma^Y_n\circ\hxD^n(f).
 \end{align*}
The conclusion follows.
\end{proof}

It follows from Lemma~\ref{L:sigD(f)} that the diagram of
Figure~\ref{Fig:FctortD} is commutative, for any $f\colon X\to Y$
in $\cI$.

\begin{figure}[htb]
 \[
 {
 \def\labelstyle{\displaystyle}
 \xymatrix{
 \hxD(Y)\ar[r]^{\sigma^Y_1} & \hxD^2(Y)\ar[r]^{\sigma^Y_2} &
 \hxD^3(Y)\ar[r]^{\sigma^Y_3} & \cdots\ \cdots\\
 \hxD(X)\ar[u]^{\hxD(f)}\ar[r]_{\sigma^X_1} &
 \hxD^2(X)\ar[u]^{\hxD^2(f)}\ar[r]_{\sigma^X_2} &
 \hxD^3(X)\ar[u]^{\hxD^3(f)}\ar[r]_{\sigma^X_3} & \cdots\ \cdots
 }
 }
 \]
\caption{A subdiagram of the functor $\tbbD$.}
\label{Fig:FctortD}
\end{figure}

We construct the functor $\tbbD\colon\cI\times\go\to\cB$ so
that all the diagrams represented by Figure~\ref{Fig:FctortD} are
subdiagrams of $\tbbD$. More precisely:
\begin{itemize}
\item[---] $\tbbD(X,n)=\hxD^{n+1}(X)$, for all
$\seq{X,n}\in\Ob\cI\times\go$.

\item[---] If $f\colon X\to Y$ in $\cI$ and $m\leq n<\omega$, let
 \[
 \tbbD(f,(m\to n))=
 \sigma^Y_n\circ\cdots\circ\sigma^Y_{m+1}\circ\hxD^{m+1}(f).
 \]
\end{itemize}
It follows from Lemma~\ref{L:sigD(f)} (that is, the commutativity
of the diagram represented on Figure~\ref{Fig:FctortD}) that we
also have
 \[
 \tbbD(f,(m\to n))=
 \sigma^Y_n\circ\cdots\circ\sigma^Y_{i+1}\circ\hxD^{i+1}(f)\circ
 \sigma^X_i\circ\cdots\circ\sigma^X_{m+1},
 \]
for all $i\in\set{m,m+1,\dots,n}$.

\section{The functor $\xQ^n$ and the equivalence $\zeta$}
\label{S:Qnzeta}

In this section we shall fix again categories $\cI$ and $\cB$,
the latter having nonempty finite products, together with a functor
$\bbD\colon\cI\to\Retr(\cB,\cB)$.
As in Section~\ref{S:D2tD}, let $\bbD$ be given by
$\seq{\xD,\hxD,\eps,\mu}$.

We assume that the unfolding
$\tbbD$ (see Section~\ref{S:D2tD}) has a lifting $\bbE$
with respect to~$\xF$. For $f\colon X\to Y$ in $\cI$ and
$n\in\gos$, we put
 \begin{align*}
 \xE^n(X)&=\bbE(X,n-1),\\
 s^X_n&=\bbE(\id_X,(n-1\to n)),\\
 \xE^n(f)&=\bbE(f,(n-1\to n-1)),
 \end{align*}
and we let the isomorphisms $\eta^X_n\colon\xF\xE^n(X)\to\hxD^n(X)$
witness the equivalence between $\xF\bbE$ and $\tbbD$. The various
relations between these morphisms can be read on the commutative
diagrams represented on Figure~\ref{Fig:xEeta}.

\begin{figure}[htb]
 \[
 {
 \def\labelstyle{\displaystyle}
 \xymatrix{
 \xE^n(Y)\ar[r]^{s^Y_n} & \xE^{n+1}(Y) &
 \xF\xE^n(X)\ar[r]^{\xF\xE^n(f)}\ar[d]_{\eta^X_n}^{\cong}
 & \xF\xE^n(Y)\ar[d]^{\eta^Y_n}_{\cong} &
 \xF\xE^n(X)\ar[r]^{\xF(s^X_n)}\ar[d]_{\eta^X_n}^{\cong} &
 \xF\xE^{n+1}(X)\ar[d]_{\eta^X_{n+1}}^{\cong}\\
 \xE^n(X)\ar[u]_{\xE^n(f)}\ar[r]_{s^X_n} &
 \xE^{n+1}(X)\ar[u]_{\xE^{n+1}(f)} &
 \hxD^n(X)\ar[r]_{\hxD^n(f)} & \hxD^n(Y) &
 \hxD^n(X)\ar[r]_{\sigma^X_n} & \hxD^{n+1}(X)
 }
 }
 \]
\caption{The lifting $\bbE$ and the equivalence $\eta$.}
\label{Fig:xEeta}
\end{figure}

Let $X\in\Ob\cI$ and let $n\in\gos$. Since the functor $\xF$ is
projectable and $\alpha^X_n$ is either an isomorphism or a
projection (thus so is $\alpha^X_n\circ\eta^X_n$), the triple
$\seq{\alpha^X_n\circ\eta^X_n,\xE^n(X),\hxD(X)}$ has a
projectability witness, say, $\seq{a^X_n,\zeta^X_n}$, with an
epic\linebreak
$a^X_n\colon\xE^n(X)\onto\xQ^n(X)$ and an isomorphism
$\zeta^X_n\colon\xF\xQ^n(X)\to\hxD(X)$ (see
Figure~\ref{Fig:GRaXnzeta}).

\begin{figure}[htb]
 \[
 {
 \def\labelstyle{\displaystyle}
 \xymatrix{
 \xF\xE^n(X)\ar[r]^{\xF(a^X_n)}\ar[d]_{\eta^X_n}^{\cong} &
 \xF\xQ^n(X)\ar[d]^{\zeta^X_n}_{\cong}\\
 \hxD^n(X)\ar[r]_{\alpha^X_n} & \hxD(X)
 }
 }
 \]
\caption{The projectability witness $\seq{a^X_n,\zeta^X_n}$
for $\seq{\alpha^X_n\circ\eta^X_n,\xE^n(X),\hxD(X)}$.}
\label{Fig:GRaXnzeta}
\end{figure}

\begin{lemma}\label{L:existsbar}
For all $X\in\Ob\cI$ and all $n\in\gos$, there exists a unique
morphism $\ol{s}^X_n\colon\xQ^n(X)\to\xQ^{n+1}(X)$ such that
$\ol{s}^X_n\circ a^X_n=a^X_{n+1}\circ s^X_n$ and
$\zeta^X_{n+1}\circ\xF(\ol{s}^X_n)=\rho_X\circ\zeta^X_n$.
\end{lemma}

The meaning of Lemma~\ref{L:existsbar} is illustrated on
Figure~\ref{Fig:olsXn}.

\begin{figure}[htb]
 \[
 {
 \def\labelstyle{\displaystyle}
 \xymatrix{
 \xE^n(X)\ar@{->>}[d]_{a^X_n}\ar[r]^{s^X_n} &
 \xE^{n+1}(X)\ar@{->>}[d]^{a^X_{n+1}} &
 \xF\xQ^n(X)\ar[d]_{\zeta^X_n}\ar[r]^{\xF(\ol{s}^X_n)} &
 \xF\xQ^{n+1}(X)\ar[d]^{\zeta^X_{n+1}}\\
 \xQ^n(X)\ar[r]^{\ol{s}^X_n} & \xQ^{n+1}(X) &
 \hxD(X)\ar[r]^{\rho_X} & \hxD(X)
 }
 }
 \]
\caption{The morphism $\ol{s}^X_n$.}
\label{Fig:olsXn}
\end{figure}

\begin{proof}
Put $\psi=(\zeta^X_{n+1})^{-1}\circ\rho_X\circ\zeta^X_n$, so
$\psi\colon\xF\xQ^n(X)\to\xF\xQ^{n+1}(X)$. We compute:
 \begin{align*}
 \psi\circ\xF(a^X_n)
 &=(\zeta^X_{n+1})^{-1}\circ\rho_X\circ\zeta^X_n\circ\xF(a^X_n)\\
 &=(\zeta^X_{n+1})^{-1}\circ\rho_X\circ\alpha^X_n\circ\eta^X_n
 &&(\text{see Figure~\ref{Fig:GRaXnzeta}})\\
 &=(\zeta^X_{n+1})^{-1}\circ\alpha^X_{n+1}\circ\sigma^X_n
 \circ\eta^X_n&&(\text{see Lemma~\ref{L:sigoalph}})\\
 &=\xF(a^X_{n+1})\circ(\eta^X_{n+1})^{-1}\circ\sigma^X_n
 \circ\eta^X_n&&(\text{see Figure~\ref{Fig:GRaXnzeta}})\\
 &=\xF(a^X_{n+1})\circ\xF(s^X_n)
 &&(\text{see Figure~\ref{Fig:xEeta}})\\
 &=\xF(a^X_{n+1}\circ s^X_n).
 \end{align*}
Hence, since $\seq{a^X_n,\zeta^X_n}$ is a projectability witness
for $\seq{\alpha^X_n\circ\eta^X_n,\xE^n(X),\hxD(X)}$, there exists
a unique morphism $\ol{s}^X_n\colon\xQ^n(X)\to\xQ^{n+1}(X)$ such
that $\ol{s}^X_n\circ a^X_n=a^X_{n+1}\circ s^X_n$ and
$\xF(\ol{s}^X_n)=\psi$. The latter condition is equivalent to
$\zeta^X_{n+1}\circ\xF(\ol{s}^X_n)=\rho_X\circ\zeta^X_n$.
\end{proof}

\begin{lemma}\label{L:Qn(f)}
Let $f\colon X\to Y$ in $\cI$ and let $n\in\gos$. There exists a
unique morphism $\xQ^n(f)\colon\xQ^n(X)\to\xQ^n(Y)$ such that
$\xQ^n(f)\circ a^X_n=a^Y_n\circ\xE^n(f)$ and
$\zeta^Y_n\circ\xF\xQ^n(f)=\hxD(f)\circ\zeta^X_n$.
\end{lemma}

The meaning of Lemma~\ref{L:Qn(f)} is illustrated on
Figure~\ref{Fig:Qnf}.

\begin{figure}[htb]
 \[
 {
 \def\labelstyle{\displaystyle}
 \xymatrix{
 \xE^n(X)\ar@{->>}[d]_{a^X_n}\ar[r]^{\xE^n(f)} &
 \xE^n(Y)\ar@{->>}[d]^{a^Y_n} &
 \xF\xQ^n(X)\ar[d]_{\zeta^X_n}\ar[rr]^{\xF\xQ^n(f)} & &
 \xF\xQ^n(Y)\ar[d]^{\zeta^Y_n}\\
 \xQ^n(X)\ar[r]^{\xQ^n(f)} & \xQ^n(Y) & \hxD(X)\ar[rr]^{\hxD(f)}
 & & \hxD(Y)
 }
 }
 \]
\caption{The morphism $\xQ^n(f)$.}
\label{Fig:Qnf}
\end{figure}

\begin{proof}
Put $\psi=(\zeta^Y_n)^{-1}\circ\hxD(f)\circ\zeta^X_n$, so
$\psi\colon\xF\xQ^n(X)\to\xF\xQ^n(Y)$. We compute:
 \begin{align*}
 \psi\circ\xF(a^X_n)&=
 (\zeta^Y_n)^{-1}\circ\hxD(f)\circ\zeta^X_n\circ\xF(a^X_n)\\
 &=(\zeta^Y_n)^{-1}\circ\hxD(f)\circ\alpha^X_n\circ\eta^X_n
 &&(\text{see Figure~\ref{Fig:GRaXnzeta}})\\
 &=(\zeta^Y_n)^{-1}\circ\alpha^Y_n\circ\hxD^n(f)\circ\eta^X_n
 &&(\text{see Lemma~\ref{L:alphD(f)}})\\
 &=\xF(a^Y_n)\circ(\eta^Y_n)^{-1}\circ\hxD^n(f)\circ\eta^X_n
 &&(\text{see Figure~\ref{Fig:GRaXnzeta}})\\
 &=\xF(a^Y_n)\circ\xF\xE^n(f)
 &&(\text{see Figure~\ref{Fig:xEeta}})\\
 &=\xF\bigl(a^Y_n\circ\xE^n(f)\bigr).
 \end{align*}
Hence, since $\seq{a^X_n,\zeta^X_n}$ is a projectability witness
for $\seq{\alpha^X_n\circ\eta^X_n,\xE^n(X),\hxD(X)}$, there exists
a unique morphism $\xQ^n(f)\colon\xQ^n(X)\to\xQ^n(Y)$ such that
$\xQ^n(f)\circ a^X_n=a^Y_n\circ\xE^n(f)$ and $\xF\xQ^n(f)=\psi$.
The latter condition is equivalent to
$\zeta^Y_n\circ\xF\xQ^n(f)=\hxD(f)\circ\zeta^X_n$.
\end{proof}

\begin{lemma}\label{L:Qnfunct}
For any $n\in\gos$, the correspondences $\xQ^n$ \pup{on objects
and morphisms} define a functor from $\cI$ to $\cA$. Furthermore,
for any $f\colon X\to Y$ in $\cI$, the equality
$\xQ^{n+1}\circ\ol{s}^X_n=\ol{s}^Y_n\circ\xQ^n(f)$ holds \pup{see
Figure~\textup{\ref{Fig:Qnf2}}}.
\end{lemma}

Lemma~\ref{L:Qnfunct} means that the correspondences
$\seq{X,n}\mapsto\xQ^{n+1}(X)$ and $\seq{f,n}\mapsto\xQ^{n+1}(f)$
define a functor from $\cI\times\go$ to $\cA$.

\begin{figure}[htb]
 \[
 {
 \def\labelstyle{\displaystyle}
 \xymatrix{
 \xQ^n(Y)\ar[r]^{\ol{s}^Y_n} & \xQ^{n+1}(Y)\\
 \xQ^n(X)\ar[u]^{\xQ^n(f)}\ar[r]^{\ol{s}^X_n}
 &\xQ^{n+1}(X)\ar[u]_{\xQ^{n+1}(f)}
 }
 }
 \]
\caption{Extending $\xQ$ to a functor from $\cI\times\go$ to
$\cA$.}
\label{Fig:Qnf2}
\end{figure}

\begin{proof}
We verify that $\xQ^n(g\circ f)=\xQ^n(g)\circ\xQ^n(f)$, for all
$f\colon X\to Y$ and $g\colon Y\to Z$ in $\cI$. We compute:
 \begin{align*}
 \xQ^n(g)\circ\xQ^n(f)\circ a^X_n&=
 \xQ^n(g)\circ a^Y_n\circ\xE^n(f)
 &&(\text{see Lemma~\ref{L:Qn(f)}})\\
 &=a^Z_n\circ\xE^n(g)\circ\xE^n(f)
 &&(\text{see Lemma~\ref{L:Qn(f)}})\\
 &=a^Z_n\circ\xE^n(g\circ f)
 &&(\text{because }\xE^n\text{ is a functor})\\
 &=\xQ^n(g\circ f)\circ a^X_n.
 \end{align*}
Since $a^X_n$ is an epic, we obtain that
$\xQ^n(g)\circ\xQ^n(f)=\xQ^n(g\circ f)$. It is clear that $\xQ^n$
preserves identities, hence it is a functor. We now prove the
additional equality:
 \begin{align*}
 \xQ^{n+1}(f)\circ\ol{s}^X_n\circ a^X_n
 &=\xQ^{n+1}(f)\circ a^X_{n+1}\circ s^X_n
 &&(\text{see Lemma~\ref{L:existsbar}})\\
 &=a^Y_{n+1}\circ\xE^{n+1}(f)\circ s^X_n
 &&(\text{see Lemma~\ref{L:Qn(f)}})\\
 &=a^Y_{n+1}\circ s^Y_n\circ\xE^n(f)
 &&(\text{see Figure~\ref{Fig:xEeta}})\\
 &=\ol{s}^Y_n\circ a^Y_n\circ\xE^n(f)
 &&(\text{see Lemma~\ref{L:existsbar}}))\\
 &=\ol{s}^Y_n\circ\xQ^n(f)\circ a^X_n.
 &&(\text{see Lemma~\ref{L:Qn(f)}})
 \end{align*}
Therefore, since $a^X_n$ is an epic,
$\xQ^{n+1}(f)\circ\ol{s}^X_n=\ol{s}^Y_n\circ\xQ^n(f)$.
\end{proof}

\section{Constructing a lifting of $\xD$}\label{S:LiftxD}

In the present section we shall make, until the statement of
Theorem~\ref{T:Main}, the same assumptions as in
Section~\ref{S:Qnzeta}. Furthermore, we assume that $\cA$ has
directed $\go$-colimits and that $\xF$ preserves
directed $\go$-colimits (see Section~\ref{S:Basic}). We keep
the notations of Section~\ref{S:Qnzeta}, in particular for $\xQ^n$
and $\zeta^X_n$.

For any $X\in\Ob\cI$, let $\xR(X)$ denote the colimit of the
sequence
 \begin{equation}\label{Eq:DiagQnR}
 {
 \def\labelstyle{\displaystyle}
 \xymatrix{
 \xQ^1(X)\ar[r]^{\ol{s}^X_1} & \xQ^2(X)\ar[r]^{\ol{s}^X_2} &
 \xQ^3(X)\ar[r]^(.4){\ol{s}^X_3} & \cdots\ \cdots,
 }
 }
 \end{equation}
with limiting morphisms $t^X_n\colon\xQ^n(X)\to\xR(X)$. Hence
$t^X_n=t^X_{n+1}\circ\ol{s}^X_n$, for all $n\in\gos$ (see
Figure~\ref{Fig:delta}).

Furthermore, it follows from Lemma~\ref{L:DirLimIdemp} that
$\xD(X)$ is the colimit of the sequence
 \[
 {
 \def\labelstyle{\displaystyle}
 \xymatrix{
 \hxD(X)\ar[r]^{\rho_X} & \hxD(X)\ar[r]^{\rho_X} &
 \hxD(X)\ar[r]^(.4){\rho_X} & \cdots\ \cdots
 }
 }
 \]
with constant limiting morphism $\mu_X\colon\hxD(X)\to\xD(X)$.
Now we apply $\xF$ to the diagram \eqref{Eq:DiagQnR}, together
with its colimit and limiting morphisms. Since $\xF$
preserves directed $\go$-colimits and by Lemma~\ref{L:existsbar}
(see Figure~\ref{Fig:olsXn}), there exists a unique isomorphism
$\delta^X\colon\xF\xR(X)\to\xD(X)$ such that
$\delta^X\circ\xF(t^X_n)=\mu_X\circ\zeta^X_n$, for all $n\in\gos$
(see Figure~\ref{Fig:delta}).

\begin{figure}[htb]
 \[
 {
 \def\labelstyle{\displaystyle}
 \xymatrix{
 \xQ^n(X)\ar[r]^{\ol{s}^X_n}\ar[dr]_{t^X_n} &
 \xQ^{n+1}(X)\ar[d]^{t^X_{n+1}} &
 \xF\xQ^n(X)\ar[rr]^{\xF(t^X_n)}\ar[d]_{\zeta^X_n} & &
 \xF\xR(X)\ar[d]^{\delta^X}\\
 & \xR(X) &\hxD(X)\ar[rr]^{\mu_X} & & \xD(X)
 }
 }
 \]
\caption{Defining $\xR(X)$, $t^X_n$, and $\delta^X$.}
\label{Fig:delta}
\end{figure}

We shall prove that $\xR$ is a lifting of $\xD$, with category
equivalence $\delta$.

First, let $f\colon X\to Y$ in $\cI$. By the universal property of
the colimit and by Lemma~\ref{L:Qnfunct} (see
Figure~\ref{Fig:Qnf2}), there exists a unique morphism
$\xR(f)\colon\xR(X)\to\xR(Y)$ such that $\xR(f)\circ
t^X_n=t^Y_n\circ\xQ^n(f)$ for all $n\in\gos$ (see
Figure~\ref{Fig:defR(f)}).

\begin{figure}[htb]
 \[
 {
 \def\labelstyle{\displaystyle}
 \xymatrix{
 \xQ^n(X)\ar[r]^{t^Y_n} & \xR(Y)\\
 \xQ^n(X)\ar[u]^{\xQ^n(f)}\ar[r]^{t^X_n} & \xR(X)\ar[u]_{\xR(f)}
 }
 }
 \]
\caption{Defining $\xR(f)$.}
\label{Fig:defR(f)}
\end{figure}

\begin{lemma}\label{L:Rfunctor}
The correspondence $X\mapsto\xR(X)$, $f\mapsto\xR(f)$ defines a
functor from~$\cI$ to~$\cA$.
\end{lemma}

\begin{proof}
As $\xR$ obviously preserves identities, it is sufficient to
verify that $\xR(g\circ f)=\xR(g)\circ\xR(f)$ whenever
$f\colon X\to Y$ and $g\colon Y\to Z$ in $\cI$. For all
$n\in\gos$, we compute, by using the definitions of $\xR(f)$ and
$\xR(g)$:
 \begin{align*}
 \xR(g)\circ\xR(f)\circ t^X_n&=\xR(g)\circ t^Y_n\circ\xQ^n(f)\\
 &=t^Z_n\circ\xQ^n(g)\circ\xQ^n(f)\\
 &=t^Z_n\circ\xQ^n(g\circ f)
 &&(\text{see Lemma~\ref{L:Qnfunct}}),
 \end{align*}
whence, by the definition of $\xR(g\circ f)$, we obtain the
desired equality.
\end{proof}

\begin{lemma}\label{L:deltaequiv}
The correspondence $X\mapsto\delta^X$ defines a category
equivalence from the composition $\xF\xR$ to~$\xD$.
\end{lemma}

In particular, $\xR$ is a lifting of $\xD$ with respect to~$\xF$.

\begin{proof}
We need to prove that the equality
$\delta^Y\circ\xF\xR(f)=\xD(f)\circ\delta^X$ holds, for any
$f\colon X\to Y$ in $\cI$. Put $\varphi=\delta^Y\circ\xF\xR(f)$
and $\psi=\xD(f)\circ\delta^X$. For all $n\in\gos$, we compute:
 \begin{align*}
 \varphi\circ\xF(t^X_n)
 &=\delta^Y\circ\xF(t^Y_n)\circ\xF\xQ^n(f)
 &&(\text{see Figure~\ref{Fig:defR(f)}})\\
 &=\mu_Y\circ\zeta^Y_n\circ\xF\xQ^n(f)
 &&(\text{see Figure~\ref{Fig:delta}})\\
 &=\mu_Y\circ\hxD(f)\circ\zeta^X_n
 &&(\text{see Lemma~\ref{L:Qn(f)}})\\
 &=\xD(f)\circ\mu_X\circ\zeta^X_n
 &&(\text{see Figure~\ref{Fig:FunctbbD}})\\
 &=\psi\circ\xF(t^X_n)
 &&(\text{see Figure~\ref{Fig:delta}}).
 \end{align*}
As $\xF\xR(X)=\varinjlim_n\xF\xQ^n(X)$ with limiting morphisms
$\xF(t^X_n)\colon\xF\xQ^n(X)\to\xF\xR(X)$, it follows from the
uniqueness statement in the definition of the colimit that
$\varphi=\nobreak\psi$.
\end{proof}

Hence we have completed the proof of our main technical result.

\begin{theorem}\label{T:Main}
Let $\cA$, $\cB$, and $\cI$ be categories, together with functors
$\xF\colon\cA\to\cB$ and $\bbD\colon\cI\to\Retr(\cB,\cB)$. We
denote by $\Pi$ the projection functor from
$\Retr(\cB,\cB)$ to~$\cB$, and we assume the following:
\begin{enumerate}
\item The category $\cA$ has directed $\go$-colimits.

\item The category $\cB$ has nonempty finite products.

\item The functor $\xF$ preserves directed $\go$-colimits.

\item The functor $\xF$ is projectable.
\end{enumerate}

If the unfolding $\tbbD$ of $\bbD$ \pup{see
Section~\textup{\ref{S:D2tD}}} has a lifting with respect to
$\xF$, then so does the composition $\Pi\bbD$.
\end{theorem}

Observe that with the notation of Section~\ref{S:D2tD},
$\Pi\bbD=\xD$.

We will always use Theorem~\ref{T:Main} in the following way. For
functors $\xF\colon\cA\to\cB$ and $\xD\colon\cI\to\cB$, we first
extend $\xD$ to a functor $\bbD\colon\cI\to\Retr(\cB,\cB)$, in the
sense that $\Pi\bbD=\xD$. This can sometimes be done by using
results in \cite{Ultra} (as in the proof of
Theorem~\ref{T:LiftConVar}(ii) in the present paper). Then, for the
associated functor $\hxD\colon\cI\to\cB$, the vertices of the
unfolding $\tbbD$ are finite powers of the form
$\hxD^n(X)$, where $X\in\Ob\cI$ and $n\in\gos$. Hence
Theorem~\ref{T:Main} reduces lifting of diagrams with vertices in
$\cB$ to lifting of diagrams whose vertices are finite powers of
objects in the range of $\hxD$. In the context of
Theorem~\ref{T:LiftConVar}(ii), the former are finite distributive
\jzs s, while the latter are finite Boolean \jzs s.

\section{Lifting diagrams of \jzs s with respect to the compact
congruence functor on a variety}\label{S:LiftCon}

For unexplained definitions of universal algebra we refer the
reader to G. Gr\"at\-zer~\cite{GrUA} or S. Burris and H.\,P.
Sankappanavar~\cite{BuSa}.

Let $\cV$ be a variety of algebras in a similarity type $\Sigma$.
So $\cV$, endowed with its homomorphisms, becomes a category. This
category has many nice properties, in particular, it has all small
limits and small colimits. For a member $A$ of $\cV$, let $\Conc A$
denote the semilattice of compact (i.e., finitely generated)
congruences of
$A$, endowed with join. Then $\Conc A$ is a \jzs. Furthermore, if
$f\colon A\to B$ is a homomorphism in $\cV$, one can define a
\jzh\ $\Conc f\colon\Conc A\to\Conc B$ by
 \[
 (\Conc f)(\la)=\text{congruence of }B\text{ generated by }
 \setm{\seq{f(x),f(y)}}{\seq{x,y}\in\la},
 \]
for any $\la\in\Conc A$. These rules determine a functor, still
denoted by $\Conc$, from~$\cV$ to the category~$\cS$ of
\jzs s and \jzh s. This functor preserves directed colimits.

The key point that relates this context to our general categorical
framework is the following.

\begin{lemma}\label{L:IdIndProj}
Let $A$ be an algebra in $\cV$, let $S$ be a \jzs, and let
$\varphi\colon\Conc A\onto S$ be an ideal-induced \jzh\ \pup{see
Section~\textup{\ref{S:Basic}}}. Then the triple
$\seq{\varphi,A,S}$ has a projectability witness with respect to
the functor $\Conc$. In particular, the functor $\Conc$ is
projectable.
\end{lemma}

\begin{proof}
Put $I=\varphi^{-1}\set{0}$ and denote by $\la$ the congruence of
$A$ defined by
 \[
 \seq{x,y}\in\la\ \Longleftrightarrow\ \varphi\Theta_A(x,y)=0,
 \text{ for all }x,\,y\in A,
 \]
where $\Theta_A(x,y)$ denotes the congruence of $A$ generated by
the pair $\seq{x,y}$. Put $\oll{A}=A/{\la}$ and denote by
$p\colon A\onto\oll{A}$ the canonical projection. Since $p$ is
surjective, it is an epic.

\begin{sclaim}
For all $\lu\in\Conc A$, $\varphi(\lu)=0$ if{f} $\lu\subseteq\la$.
\end{sclaim}

\begin{scproof}
Suppose first that $\varphi(\lu)=0$, and let $x$, $y\in A$ such
that $x\equiv y\pmod{\lu}$. This can be written
$\Theta_A(x,y)\subseteq\lu$, thus, applying $\varphi$ and the
assumption, $\varphi\Theta_A(x,y)=0$, that is,
$x\equiv y\pmod{\la}$. Hence $\lu\subseteq\la$.
Conversely, suppose that $\lu\subseteq\la$. Since $\lu$ is
compact, it can be written in the form
$\lu=\bigvee_{i<n}\Theta_A(x_i,y_i)$, where $n<\go$ and $x_i$,
$y_i\in A$ for all $i<n$. For all $i<n$, the relation
$x_i\equiv y_i\pmod{\la}$ holds, so $\varphi\Theta_A(x_i,y_i)=0$.
Joining over $i<n$ yields $\varphi(\lu)=0$.
\end{scproof}

Now let $\lx$, $\ly\in\Conc A$. Since $\varphi$ is ideal-induced,
the following equivalences hold:
 \begin{align*}
 (\Conc p)(\lx)=(\Conc p)(\ly)&\Leftrightarrow
 \lx\vee\la=\ly\vee\la
 &&(\text{see \cite[Theorem~I.11.3]{GrUA}})\\
 &\Leftrightarrow\varphi(\lx)=\varphi(\ly)
 &&(\text{by the Claim above}).
 \end{align*}
In particular, there exists a unique isomorphism
$\eps\colon\Conc\oll{A}\to S$ such that $\varphi=\eps\circ\Conc p$.

Now let $f\colon A\to B$ a homomorphism in $\cV$ and let
$\eta\colon\Conc\oll{A}\to\Conc B$ such that
$\Conc f=\eta\circ\Conc p$. For all $x$, $y\in A$, if
$x\equiv y\pmod{\la}$, then $p(x)=p(y)$, thus
$\Theta_B(f(x),f(y))=
(\Conc f)\Theta_A(x,y)=\eta\Theta_A(p(x),p(y))=0$, so $f(x)=f(y)$.
Hence $\ker p\subseteq\ker f$, and thus there exists a unique
homomorphism $g\colon\oll{A}\to B$ such that $f=g\circ p$. Thus
$\eta\circ(\Conc p)=\Conc f=(\Conc g)\circ(\Conc p)$, and hence,
since $\Conc p$ is surjective, $\eta=\Conc g$. Therefore,
$\seq{p,\eps}$ is a projectability witness for $\seq{\varphi,A,S}$.

For \jzs s $S_0$ and $S_1$, the canonical projection from
$S_0\times S_1$ onto~$S_0$ is obviously ideal-induced. Hence, by
the result above, the functor $\Conc$ is projectable.
\end{proof}

We shall now present some applications of
Lemma~\ref{L:IdIndProj}. These applications will depend on
a few results from~\cite{Ultra}.
We first introduce the following subcategories of~$\cS$:
\begin{itemize}
\item $\cM$, the subcategory of all \jzs s and \jze s.

\item $\cS_{\fin}$, the full subcategory of all finite members of
$\cS$.

\item $\cS^*_{\fin}$, the full subcategory of all finite,
simple, atomistic lattices and \jzh s.

\item $\cD$, the full subcategory of all distributive \jzs s.

\item $\cB$, the full subcategory of all Boolean \jzs s.
\end{itemize}

Now we quote the corresponding results obtained in~\cite{Ultra}.
While Proposition~\ref{P:GSretr} is obtained, in~\cite{Ultra}, by a
simple application of a construction by G. Gr\"atzer and E.\,T.
Schmidt in~\cite{GrSc99}, we know no proof of
Proposition~\ref{P:FRSemil} substantially simpler than the
relatively involved categorical argument presented in
\cite{Ultra}. The importance of finite simple atomistic lattices
(or, more generally, finite simple lattices whose atoms join to
the unit) for congruence representation problems of algebras is
illustrated by P.\,P. P\'alfy and P. Pudl\'ak~\cite{P5}.

\begin{proposition}\label{P:GSretr}
There exists a functorial retraction $\seq{\Upsilon,\xi,\rho}$ of
the category $\cS_{\fin}\cap\nobreak\cM$ to the full subcategory
$\cS^*_{\fin}\cap\cM$. Furthermore, $\xi_K$ is a \ble, for every
finite \jzs\ $K$.
\end{proposition}

\begin{proposition}\label{P:FRSemil}
There exists a functorial retraction $\seq{\Phi,\eps,\mu}$ of the
category $\cD_{\fin}\cap\nobreak\cM$ to the category
$\cB_{\fin}\cap\cM$. Furthermore, $\eps_A$ is a \jzue, for all
$A\in\Ob\cD$.
\end{proposition}

Now the promised applications.

\begin{theorem}\label{T:LiftConVar}
Let $\cV$ be a variety of algebras. Then the following statements
hold:
\begin{enumerate}
\item If every diagram of finite products of
finite atomistic simple lattices and \jzue s, indexed by a
lattice, can be lifted with respect to the $\Conc$ functor on
$\cV$, then so can every diagram of finite \jzus s and
\jzue s, indexed by a lattice.

\item If every diagram of finite Boolean \jzs s and
\jze s, indexed by a lattice, can be
lifted with respect to the $\Conc$ functor on $\cV$, then so can
every diagram of finite distributive \jzs s and \jze s, indexed by
a lattice.

\item The analogue of \textup{(ii)} above for
\jzue s also holds.
\end{enumerate}
\end{theorem}

We shall discuss later why we do not formulate the analogue of
Theorem~\ref{T:LiftConVar}(i) for \jze s (although it holds!), see
Proposition~\ref{P:NonLiftSig}.

\begin{proof}
We provide a proof for (i); it uses
Proposition~\ref{P:GSretr}. The proofs of (ii) and (iii) are
similar. Of course, the proof of (ii) uses
Proposition~\ref{P:FRSemil} and the triple $\seq{\Phi,\eps,\mu}$
instead of Proposition~\ref{P:GSretr} and the triple
$\seq{\Upsilon,\xi,\rho}$.

Let $F\colon\cV\to\cS$ be the functor that with an algebra $A$
associates $\Conc A$ and that with a homomorphism $f$ associates
$\Conc f$. Let $\cI$ be a lattice, viewed as a
category, and let $\xD\colon\cI\to\cS_{\fin}\cap\cM$ be a
$\cI$-indexed diagram of finite \jzus s and \jzue s. We define a
functor $\bbD\colon\cI\to\Retr(\cS,\cS)$ by putting
 \begin{align*}
 \bbD(x)&=\seq{\xD(x),\Upsilon\xD(x),\xi_{\xD(x)},\rho_{\xD(x)}},\\
 \bbD(x\to y)&=\seq{\xD(x\to y),\Upsilon\xD(x\to y)}
 \end{align*}
for all $x\leq y$ in $\cI$. Observe that $\xD(x)\in\Ob\cS_{\fin}$,
$\Upsilon\xD(x)\in\Ob\cS^*_{\fin}$, and both $\xD(x\to y)$ and
$\Upsilon\xD(x\to y)$ are \jzue s.

Now we take a closer look at the unfolding
$\tbbD\colon\cI\times\go\to\cS$ of $\bbD$.
Its objects are finite products of finite simple lattices.
Its arrows are finite compositions of arrows
of the form either $\sigma^x_n$ (for $x\in\cI$ and $n\in\gos$),
which is a \emph{section} of $\cS$ (for
$\pi^x_{n+1}\circ\sigma^x_n$ is an identity), or
$\hxD^n(x\to y)$, for $x\leq y$ in $\cI$ and $n\in\gos$. As
$\hxD(x\to y)=\Upsilon\xD(x\to y)$ is a \jzue\ and $f\vpd g$ is a
\jzue\ whenever both $f$ and $g$ are \jzue s, $\hxD^n(x\to y)$ is a
\jzue\ for all $n$. Consequently, all arrows of $\tbbD$ are
\jzue s.

Hence, by assumption, the diagram $\tbbD$ has a lifting with
respect to the $\Conc$ functor. By Lemma~\ref{L:IdIndProj}, the
$\Conc$ functor is projectable. Moreover, $\cV$ has arbitrary
directed colimits, $\Conc$ preserves all directed colimits, and
$\cS$ has finite products. Therefore, by Theorem~\ref{T:Main}, the
diagram $\xD=\Pi\bbD$ has a lifting with respect to $\Conc$.
\end{proof}

The statement of Theorem~\ref{T:LiftConVar} can be modified in
many ways, especially on the category $\cI$ indexing the diagram
$\xD$. Indeed, the unfolding $\tbbD$ is indexed by
$\cI\times\omega$. So, instead of requiring that $\cI$ be a
lattice, we could have required that $\cI$ be a poset, or a
distributive lattice, or even an arbitrary category. The
only restriction on the arrows of $\xD$ is that they are required
to be \emph{\jze s} (\jzue s in (i)), in order to be able to apply
Propositions~\ref{P:GSretr} and~\ref{P:FRSemil}.

On the other hand, although the corresponding analogues of
Theorem~\ref{T:LiftConVar} are still valid, they may be so
vacuously, that is, their assumptions may not hold.

An example of the latter situation is the analogue of
Theorem~\ref{T:LiftConVar}(i) for \jze s (not necessarily
preserving the unit). Let us be more specific.

\begin{example}\label{Ex:LampLim}
Let $\cV$ be the variety of all algebras with given similarity
type~$\Sigma$. It follows from results in
R. Freese, W.\,A. Lampe, and W. Taylor~\cite{FLT} that there
exists a \jzs\ $S$ that is not isomorphic to $\Conc A$, for any
member~$A$ of~$\cV$. In fact, $S$ may be defined as the \jzs\ of
subspaces of an infinite-dimensional vector space over a field
with large enough cardinality. On the other hand, $S$ is the
directed union of its finite \jz-subsemilattices. Since the
$\Conc$ functor preserves arbitrary directed colimits, the
corresponding diagram of \jzs s has no lifting with respect to
$\Conc$. This diagram may be taken indexed by the poset of all
finite subsets of~$S$, which is a distributive lattice. Therefore,
by applying the analogue of Theorem~\ref{T:LiftConVar}(i) for
\jze s, we obtain the following result.
\end{example}

\begin{proposition}\label{P:NonLiftSig}
Let $\Sigma$ be a similarity type of algebras. Then there exists a
diagram, indexed by a distributive lattice, of finite
products of finite simple lattices and \jze s, without a lifting,
with respect to the $\Conc$ functor, by algebras with similarity
type $\Sigma$.
\end{proposition}

Note that there is a very simple diagram (indexed by a
finite category) of finite Boolean \jzus s and \jzue s that cannot
be lifted, with respect to the $\Conc$ functor, by algebras in any
given similarity type (see J. T\r{u}ma and F.
Wehrung, \cite[Theorem~8.1]{TuWe1} or~\cite[Introduction]{Bowtie}).
Furthermore, it is proved in J. T\r{u}ma and F. Wehrung
\cite{Bowtie} that there exists a diagram, indexed by a finite
poset, of finite Boolean \jzus s with \jzue s, that cannot be
lifted, with respect to the $\Conc$ functor, by lattices and
lattice homomorphisms. This is important in telling us how far we
may go in the formulation of the coming
Proposition~\ref{P:UltrabRed}.

\begin{example}\label{Ex:CLP}
It is still an open problem, dating back to 1945 and usually
referred to as CLP, whether every distributive \jzs\ is isomorphic
to $\Conc L$ for some lattice $L$. A survey of some of the latest
attempts on that problem can be found in~\cite{CLPSurv}. All
these attacks aim at lifting not only
\emph{objects} (distributive \jzs s), but \emph{diagrams} (of
finite distributive \jzs s). Some partial results and many
calculations suggest that it is less difficult to lift diagrams of
finite Boolean \jzs s and \jze s than to lift diagrams of
arbitrary distributive \jzs s (this is also illustrated in J.
T\r{u}ma~\cite{Tuma93}). As every distributive \jzs\ is the
directed union of its finite \jz-subsemilattices (see Fact~4,
page~100 in P. Pudl\'ak~\cite{Pudl85}), a direct application of
Theorem~\ref{T:LiftConVar}(ii) yields the following.
\end{example}

\begin{proposition}\label{P:UltrabRed}
If every diagram, indexed by a distributive
lattice, of finite Boolean \jzs s \pup{resp., \jzus s} 
with \jze s \pup{resp., \jzue s} can be
lifted, with respect to the $\Conc$ functor, by lattices and
lattice homomorphisms, then every distributive \jzs\ \pup{resp.,
\jzus} is isomorphic to $\Conc L$ for some lattice $L$.
\end{proposition}

Various other formulations are possible, for example with \jzs s
and $0$-lattice homomorphisms, or bounded \jzus s and bounded
lattices, with similar proofs. However, it is an open problem
whether any of the corresponding assumptions actually holds.

\section{Nonstable K-theory of von~Neumann regular rings}
\label{S:V(R)}

For a (unital, associative) ring $R$, let $\FP(R)$ denote the class
of all finitely generated projective right $R$-modules. For
$X\in\FP(R)$, denote by $[X]$ the isomorphism class of $X$.
Isomorphisms classes can be added, with the rule
$[X]+[Y]=[X\oplus Y]$, and the monoid $V(R)$ of all isomorphism
classes of members of $\FP(R)$ is a conical \cm\ with order-unit.

Now we specialize our class of rings. A ring $R$
is (von~Neumann) \emph{regular}, if it satisfies the axiom
$(\forall x)(\exists y)(xyx=x)$. For a regular ring $R$, the monoid
$V(R)$ is a \emph{refinement monoid}, see K.\,R.
Goodearl~\cite[Theorem~2.8]{Gvnrr}. However, there exists a
conical refinement monoid (and even the positive cone of a
dimension group) that is not isomorphic to $V(R)$, for any regular
ring $R$, see F. Wehrung~\cite{Wehr98}. This counterexample has
cardinality $\aleph_2$. The corresponding problem for smaller
cardinalities is still open, and probably very difficult, even in
the countable case, see the fundamental open problem
in K.\,R. Goodearl~\cite{Good94}. We formulate the corresponding
question for the countable case: \emph{Is every countable conical
refinement monoid with order-unit isomorphic to $V(R)$ for some
regular ring~$R$?}

In the present section we wish to provide a new insight on this
problem, based on Theorem~\ref{T:Main}. First of all, observe that
for regular rings $R$ and $S$, every ring homomorphism
$f\colon R\to S$ induces a unique monoid homomorphism
$V(f)\colon V(R)\to V(S)$ such that $V(f)([xR])=[f(x)S]$, for all
$x\in R$. Hence $V$ is a \emph{functor} from the category $\cR$ of
regular rings and (unital) ring homomorphisms to the category
$\cM$ of conical refinement monoids with order-unit and
unit-preserving homomorphisms. It is well-known that $\cR$ has
arbitrary directed colimits and that the $V$ functor preserves all
directed colimits.

We need the following lemma.

\begin{lemma}\label{L:VfunctProj}
Let $R$ be a regular ring, let $M$ be a conical refinement monoid,
and let $\varphi\colon V(R)\onto M$ be an ideal-induced monoid
homomorphism. Then the triple $\seq{\varphi,R,M}$ has a
projectability witness with respect to the functor $V$. In
particular, the functor~$V$ is projectable.
\end{lemma}

\begin{proof}
Put $J=\varphi^{-1}\set{0}$ and $I=\setm{x\in R}{[xR]\in J}$. So
$I$ is a two-sided ideal of $R$ and $J=\setm{[xR]}{x\in I}$, see
Corollary~4.4 and Proposition~4.6 in F. Wehrung~\cite{Wehr99}.
Put $\ol{R}=R/I$ and denote by $p\colon R\onto R/I$ the canonical
projection. Since $p$ is surjective, it is an epic. Denote by
$\pi\colon V(R)\onto V(R)/J$ the canonical projection.
By describing explicitly the isomorphism $V(R/I)\cong V(R)/V(I)$
(see Proposici\'o~4.1.6 and Teorema~4.1.7 in
E. Pardo~\cite{PardTh}, or P. Ara
\emph{et al.}~\cite[Proposition~1.4]{AGPO}), we obtain that there
exists an isomorphism $\zeta\colon V(\ol{R})\to V(R)/J$ such that
 \[
 \zeta\bigl([(x+I)\ol{R}]\bigr)=\pi\bigl([xR]\bigr),
 \text{ for all }x\in R.
 \]
Hence $\zeta\circ V(p)=\pi$. Since $\varphi$ is ideal-induced, it
induces a monoid isomorphism $\psi\colon V(R)/J\to M$ such
that $\psi\circ\pi=\varphi$. Hence $\eps=\psi\circ\zeta$ is an
isomorphism from $V(\ol{R})$ onto $M$ such that
$\varphi=\eps\circ V(p)$.

Finally, let $S$ be a regular ring, let $f\colon R\to S$ be a
ring homomorphism, and let $\eta\colon V(\ol{R})\to V(S)$ be a
monoid homomorphism such that $V(f)=\eta\circ V(p)$. {}From the
obvious fact that $x=0$ if{f} $[xR]=0$ (in any ring), we obtain
that $\ker p$ is contained in $\ker f$, thus there exists a unique
ring homomorphism $g\colon\ol{R}\to S$ such that $g\circ p=f$.
Hence, $V(g)\circ V(p)=V(f)=\eta\circ V(p)$, and hence, since
$V(p)$ is surjective, $V(g)=\eta$. Therefore, $\seq{p,\eps}$ is a
projectability witness for $\seq{\varphi,R,M}$.

For \cm s $M_0$ and $M_1$, the canonical projection from
$M_0\times M_1$ onto $M_0$ is obviously ideal-induced. Hence, by
the result above, the functor $V$ is projectable.
\end{proof}

Applying Theorem~\ref{T:Main}, and then
Lemma~\ref{L:VfunctProj} just as Lemma~\ref{L:IdIndProj} is used in
the proof of Theorem~\ref{T:LiftConVar}, yields the following
result.

\begin{proposition}\label{P:V(R)attempt}
Let $M$ be a conical refinement monoid with order-unit and
let~$\cM'$ be a full subcategory of $\cM$ satisfying the following
conditions:
\begin{enumerate}
\item $\cM'$ is closed under finite nonempty products.

\item Every diagram, indexed by $\go\times\go$, of members of
$\cM'$ with \emph{embeddings} \pup{see
Section~\textup{\ref{S:Basic}}}, has a lifting, with respect to the
$V$ functor, by regular rings.

\item There exists a $\go$-indexed diagram in $\Retr(\cM,\cM')$ of
the form
 \[
 {
 \def\labelstyle{\displaystyle}
 \xymatrix{
 N_0\ar@{^(->}[r]^{g_0}\ar@{->>}[d]<.5ex> &
 N_1\ar@{^(->}[r]^{g_1}\ar@{->>}[d]<.5ex> &
 N_2\ar@{^(->}[r]^(.4){g_2}\ar@{->>}[d]<.5ex> & \cdots\ \cdots\\
 M_0\ar@{^(->}[r]^{f_0}\ar[u]<.5ex> &
 M_1\ar@{^(->}[r]^{f_1}\ar[u]<.5ex>&
 M_2\ar@{^(->}[r]^(.4){f_2}\ar[u]<.5ex> & \cdots\ \cdots,
 }
 }
 \]
with all $f_n$-s and $g_n$-s \emph{embeddings}, such that
$\varinjlim_nM_n\cong M$.

\end{enumerate}
Then $M$ is isomorphic to $V(R)$, for some regular ring $R$.
\end{proposition}

However, unlike the situation with finite distributive and Boolean
semilattices (see Proposition~\ref{P:UltrabRed}), we do not have
at the present time any reasonable candidate for the category
$\cM'$. A possibility would be to use colimits of diagrams
of the form
 \[
 {
 \def\labelstyle{\displaystyle}
 \xymatrix{
 F\ar@{^(->}[r]^f & F\ar@{^(->}[r]^f & F\ar@{^(->}[r]^(.4){f}
 & \cdots\ \cdots,
 }
 }
 \]
where $F$ is a finitely generated \cm\ and $f\colon M\into M$
is an embedding from $M$ into $M$
(we need to ensure that the colimit has refinement). A related
construction is used, in C. Moreira dos Santos~\cite{CMor2}, to
find a refinement monoid whose maximal antisymmetric quotient is
not a refinement monoid. Also, our insistence on
\emph{embeddings} instead of just one-to-one monoid homomorphisms
is purely experimental---in particular, the analogue of
Proposition~\ref{P:V(R)attempt} for one-to-one monoid
homomorphisms is also valid.

\section{Concluding remarks}\label{S:Pbs}

Sections~\ref{S:LiftCon} and~\ref{S:V(R)} were written with some
amount of detail, in particular proving the projectability of both
functors $\Conc$ (on algebras) and $V$ (on regular rings). There
are many other situations in which the problem of the range of a
given projectable functor is raised. Let us mention a few
more examples of projectable functors and some of the corresponding
problems:

\begin{itemize}
\item[---] The functor $\bL$ from von~Neumann regular rings to
complemented modular lattices. For each regular ring $R$, $\bL(R)$
is the lattice of all principal right ideals of $R$. Determining
the range of the $\bL$ functor is a hard open problem,
initiated by von~Neumann's Coordinatization Theorem. Recently,
the author proved that even for countable regular rings, the range
of $\bL$ is not a first-order class, see F. Wehrung~\cite{Coord}.

\item[---] The forgetful functor that with every modular
ortholattice $L$ associates the corresponding lattice. Some
partial results about the range of that functor are established in
\cite{Coord}.

\item[---] The functor $\Dim$ that with a lattice $L$ associates
its \emph{dimension monoid} $\Dim L$, see F. Wehrung~\cite{WDim}.
Determining the range of $\Dim$ seems to be an even harder problem
than CLP (see Section~\ref{S:LiftCon}). When restricted to
complemented modular lattices, it is strongly related to the $V(R)$
representation problem discussed in Section~\ref{S:V(R)}.

\item[---] The functor $\nabla$ on refinement monoids. For a
refinement monoid $M$, $\nabla(M)$ is the maximal semilattice
quotient of $M$; it is a distributive \jzs. The
problem of the determination of the range of $\nabla$ on certain
classes of refinement monoids was initiated in K.\,R.
Goodearl and F. Wehrung~\cite{GoWe1}. An important milestone is P.
R\r{u}\v{z}i\v{c}ka's result that provides a
distributive \jzus\ of cardinality $\aleph_2$ which is not
isomorphic to $\nabla(G^+)$ for any dimension group $G$, see
P. R\r{u}\v{z}i\v{c}ka~\cite{Ruzi1}. This result was later improved
by the author to the cardinality $\aleph_1$ and larger classes of
monoids, see F. Wehrung~\cite{Wehr04}.

\end{itemize}

W.\,A. Lampe proved in~\cite{Lamp82} a result that
implies that every \jzus\ is isomorphic to $\Conc G$ for some
groupoid (i.e., nonempty set with a binary operation) $G$. The
proof presented in \cite{Lamp82} is apparently not functorial. This
raises the question whether the assumption of
Theorem~\ref{T:LiftConVar}(i) holds for some variety. The boldest
that we may ask in that direction is the following.

\begin{problem}\label{Pb:FuncLampe}
Denote by $\cG$ the category of all groupoids and groupoid
homomorphisms. Does there exist a functor $\Psi$ from the category
of \jzus s with \jzue s to the category $\cG$ such that $\Conc\Psi$
is equivalent to the identity?
\end{problem}

By a suitable analogue of Theorem~\ref{T:LiftConVar}(i), it is
sufficient to solve the problem on the category of finite products
of simple lattices with \jzue s.

A related problem is the following.

\begin{problem}\label{Pb:FinAlg}
Can every finite diagram of finite \jzs s with \jze s be lifted,
with respect to the $\Conc$ functor, by a diagram of algebras
\pup{in some similarity type}?
\end{problem}

Our next problem is related to lifting problems of distributive
\emph{lattices} with zero by various functors. Some known results
are the following:
\begin{itemize}
\item Every distributive lattice $D$ with zero is isomorphic to
$\nabla(\QQ\langle D\rangle^+)$, for a certain dimension vector
space $\QQ\langle D\rangle$ (``temperate power'') that is easily
seen to depend functorially on $D$; see \cite{GoWe1}.

\item Every distributive lattice $D$ with zero is isomorphic to
$\Conc L_D$, for a lattice $L_D$ that is a directed
$\seq{\vee,\wedge,0}$-colimit of finite atomistic lattices and
that depends functorially on $D$; see \cite{Pudl85}.

\item Every distributive lattice $D$ with zero is isomorphic to
the compact (i.e., finitely generated) ideal lattice $\Idc R$ of
some locally matricial ring $R$ (over any given field); see P.
R\r{u}\v{z}i\v{c}ka~\cite{Ruzi2}.
\end{itemize}

A ring is \emph{locally matricial}, if it is a directed colimit of
a diagram of finite products of the form
$\prod_{i<n}\mathrm{M}_{m_i}(F)$, where $F$ is a given field, $n$
and the $m_i$-s are positive integers, and $\mathrm{M}_{m_i}(F)$
denotes the $m_i\times m_i$ full matrix ring over $F$. We ask
whether ``the best of all worlds'' is reachable:

\begin{problem}\label{Pb:FunctPavel}
Let $F$ be a field. Does there exist a functor $\Phi$, from
distributive $0$,$1$-lattices with
$\seq{\vee,\wedge,0,1}$-embeddings to locally matricial rings over
$F$ with \pup{unital} ring homomorphisms, such that
$\Idc\Phi$ is equivalent to the identity?
\end{problem}

In view of some results of \cite{GoWe1}, it is reasonable to ask
whether if such a functor existed, the composition $K_0\Phi$ could
be equivalent to the functor $\QQ\langle{}_{-}\rangle$.

\section*{Acknowledgment}

This work was partially completed while the author was visiting
the Charles University (Prague). Excellent conditions provided by
the Department of Algebra are greatly appreciated.

\end{document}